\documentclass[12pt]{article}
\usepackage{fancyhdr, array,calc,graphicx,url,tabularx}
\usepackage{geometry}
\usepackage{latexsym}
\usepackage{amssymb}
\usepackage{amsmath}
\usepackage{makeidx}
\usepackage{enumerate}
\usepackage{verbatim}
\usepackage{tikz}
\usepackage[colorlinks=true,linkcolor=blue]{hyperref}
\usepackage{graphicx}
\usepackage{epsfig}
\usepackage{float}

\usepackage{changepage}

\makeindex

{
   \end{minipage}
   \vspace*{\stretch{3}}
   \clearpage
}

\renewcommand{\gg}{\gamma}

\newcommand{\bR}{{\mathbb{R}}}

\newcommand{\rest}{\restriction}

\newcommand{\card}[1]{{\vert #1 \vert} }

\newcommand{\forces}{\Vdash}

\renewcommand{\models}{\vDash}
\newcommand{\powerset}{{\wp}}
\newcommand{\cp}{{\rm crit }}

\newcommand{\cf}{{\rm cf}}

\newtheorem{theorem}{Theorem}[section]
\newtheorem{proposition}[theorem]{Proposition}
\newtheorem{definition}[theorem]{Definition}

\newtheorem{lemma}[theorem]{Lemma}
\newtheorem{corollary}[theorem]{Corollary}

\newtheorem{conjecture}[theorem]{Conjecture}

\newtheorem{remark}[theorem]{Remark}

\numberwithin{figure}{section}

\newenvironment{proof}{{\it{
Proof.}}}{\nopagebreak\mbox{}{\hfill$\square$}
\par\bigskip}
\newenvironment{Claim}{\noindent{\it
Claim} {}}{}
\newenvironment{Case}{\noindent{\it
Case}}{}

\newcommand{\rcon}[1]{Conjecture~\ref{#1}}

\newcommand{\rprop}[1]{Proposition~\ref{#1}}
\newcommand{\rthm}[1]{Theorem~\ref{#1}}
\newcommand{\rlem}[1]{Lemma~\ref{#1}}

\newcommand{\rcor}[1]{Corollary~\ref{#1}}
\newcommand{\rdef}[1]{Definition~\ref{#1}}

\newcommand{\rsec}[1]{Section~\ref{#1}}

\def\inseg{\trianglelefteq}

\def\k{\kappa}
\def\a{\alpha}
\def\b{\beta}
\def\d{\delta}

\def\l{\lambda}

\def\P{{\mathcal{P} }}
\def\W{{\mathcal{W} }}
\def\Q{{\mathcal{ Q}}}

\def\K{{\mathcal{ K}}}

\def\R{{\mathcal R}}
\def\X{{\mathbb X}}
\def\H{{\rm{HOD}}}
\def\M{{\mathcal{M}}}
\def\N{{\mathcal{N}}}

\def\T {{\mathcal{T}}}
\def\U{{\mathcal{U}}}
\def\S{{\mathcal{S}}}
\def\V{{\mathcal{V}}}
\def\X{{\mathcal{X}}}
\def\Y{{\mathcal{Y}}}

\def\card#1{\left|#1\right|}

\def\iff{\mathrel{\leftrightarrow}}

\def\and{\mathrel{\kern1pt\&\kern1pt}}

\def\inseg{\triangleleft}
\def\insegeq{\trianglelefteq}

\def\<#1>{\langle\,#1\,\rangle}

 \input xy
 \xyoption{all}
       
\title{Covering with Chang models over derived models\thanks{2000 Mathematics Subject Classifications:
03E15, 03E45, 03E60.}
\thanks{Keywords: Mouse, inner model theory, descriptive set theory, hod mouse.}
\thanks{The author's research was partially supported by the NSF Career Award DMS-1352034.}}

\author{Grigor Sargsyan}

\date{\today}

\begin{document}

\maketitle

\begin{abstract}
We present a covering conjecture that we expect to be true below superstrong cardinals. We then show that the conjecture is true in hod mice. This work is a continuation of the work that started in \cite{sargsyan2013covering}, and the main conjecture of the current paper is a revision of the UB-Covering Conjecture of \cite{sargsyan2013covering}.
\end{abstract}

One of the main projects of inner model theory is to identify canonical structures like G\"odel's $\sf{L}$ that are \textit{close} to the universe of sets. While what exactly \textit{closeness} should mean is open to interpretation, perhaps the simpliest way of interpreting it is by demanding that the successor of some singular cardinal is computed correctly. Perhaps the most well-known result of this kind is Jensen's covering lemma. A weaker version of it says that if  $0^\#$ doesn't exist then $\sf{L}$ computes the successor of any singular cardinal correctly, i.e., if $\kappa$ is singular then $(\kappa^+)^{\sf{L}}=\kappa^+$. In this paper, our goal is to introduce a covering principle that is based on the idea that the information coded into the universally Baire sets is enough to describe the successor cardinals. To state our conjecture we need to introduce \textit{universally Baire sets} and \textit{extenders}.

Recall from \cite{UB} that a set of reals is \textit{universally Baire} if all of its preimages in Hausdorff spaces have the property of Baire. As shown in \cite{UB}, $A\subseteq \bR$ is universally Baire if and only if for each cardinal $\k$ there is a pair of trees $(T, S)$\footnote{A tree on $\k$ is a downward closed subset of $\cup_{n<\omega}\omega^n\times \kappa^n$. If $T$ is a tree then $[T]$ are the set of branches of $T$ and $p[T]$ is the projection of $[T]$ on the first coordinate.} on $\k$ such that $p[T]=A$ and for all $<\k$-generics $g$\footnote{This means that the poset in question has size $<\k$.}, $V[g]\models p[T]=p[S]^c$. If $A$ has the property on the right side of the equivalence then we say that $A$ is $\k$-universally Baire. Clearly, $A$ is universally Baire if and only if it is $\k$-universally Baire for all $\k$. 

 Following Woodin (see \cite{SealingTheorem}), we let $\Gamma^\infty$ be the set of universally Baire sets. Given a $V$-generic $g$ we set $\Gamma^\infty_g=(\Gamma^\infty)^{V[g]}$ and $\bR_g=\bR^{V[g]}$. Given a universally Baire set $A$ and a generic $g$, we let $A^g$ be the canonical interpretation of $A$ in $V[g]$. More precisely, fixing $\k$ larger than the size of the poset and a pair of trees $(T, S)$ witnessing that $A$ is $\k$-universally Baire, $A^g=(p[T])^{V[g]}$\footnote{Absoluteness implies that $A^g$ is independent of $\k$ and $(T, S)$. See \cite{UB}.}. 

An extender is a system of ultrafilters. More precisely, we say $E$ is a $(\k, \l)$-extender over $M$ if there is an embedding $j: M\rightarrow N$ such that 
\begin{enumerate}
\item $\cp(j)=\k$,
\item $j(\k)\geq \l$,
\item $E=\{(a, A): a\in \lambda^{<\omega}, A\in \powerset(\k^{\card{a}})\cap M$ and $a\in j(A)\}$,
\item $N=\{ j(f)(a): f\in M, f:\kappa^{\card{a}}\rightarrow M$ and $a\in \lambda^{<\omega}\}$. 
\end{enumerate}
Notice that for each $a\in \l^{<\omega}$, $E_a=\{A: (a, A)\in E\}$ is an ultrafilter measuring $\powerset(\kappa^{\card{a}})\cap M$. The extenders we have defined are usually called \textit{short}, where shortness refers to the fact that $j(\kappa)\geq \l$. The fact that  $E_a$ is an ultrafilter on $\k$ is a consequence of shortness. Long extenders have measures concentrating on more than one cardinal. Just like ultrafilters, extenders can be defined without a reference to the embedding $j$. In our set up, $N$ is the ultrapower of $M$ by $E$, and it can be shown to be completely determined by the pair $(M, E)$. 

From now on, when we say we work in the short extender region we mean that we are in the region of large cardinals that can be defined using only short extenders. Large cardinal notions such as strong cardinals, Woodin cardinals, Shelah cardinals, superstrong cardinals and subcompact cardinals are all in the short extender region, while supercompactness is not. We will drop \textit{short} from now on, and to state our results, will use Steel's $\sf{NLE}$ (see \cite{normalization_comparison}), which stands for ``no long extender". 

\begin{conjecture}[Covering with Chang Models]\label{main conjecture} Assume $\sf{NLE}$ and suppose there are unboundedly many Woodin cardinals and strong cardinals. Let $\kappa$ be a limit of Woodin cardinals and strong cardinals and such that either $\k$ is a measurable cardinal or $\cf(\k)=\omega$. Then there is a transitive model $M$ of $\sf{ZFC-Powerset}$ such that
\begin{enumerate}
\item $Ord\cap M=\k^+$,
\item $M$ has a largest cardinal $\nu$,
\item for any $g\subseteq Coll(\omega, <\k)$, letting $\mathbb{R}^* = \bigcup_{\alpha<\k} \mathbb{R}^{V[g\cap Coll(\omega,\alpha)]}$ and $\Gamma^*=\{ A^g\cap \bR^*: \exists \a<\k(A\in \Gamma^\infty_{g\cap Coll(\omega,\alpha)})\}$, in $V(\bR^*)$, 
\begin{center}
$L(M, \bigcup_{\a<\nu}\a^\omega, \Gamma^*, \bR^*)\models \sf{AD}$.
\end{center}
\item If in addition there is no inner model with a subcompact cardinal then
\begin{center}
 $L(M)\models {\sf{ZFC}}+\powerset(\nu)=\powerset(\nu)^M+\square_\nu$.
 \end{center}
\end{enumerate}
\end{conjecture}
We state the following two immediate corollaries of $\sf{Covering\ with\ Chang\ Models}$.

\begin{corollary}\label{cor 1 square} Assume $\sf{Covering\ with\ Chang\ Models}$ holds and there is no inner model with a subcompact cardinal. Assume further that there are unboundedly many Woodin cardinals and strong cardinals. Let $\kappa$ be a limit of Woodin cardinals and strong cardinals such that either $\k$ is a measurable cardinal or $\cf(\k)=\omega$. Then $\square(\k^+)$ holds. 
\end{corollary}
\begin{proof} Let $M$ be as in  $\sf{Covering\ with\ Chang\ Models}$ applied to $\kappa$ and let $\nu$ be its largest cardinal. Then as $L(M)\models \square_\nu$, we have  a sequence $\vec{C}=(C_\alpha: \a<\kappa^+)\in L(M)$ such that
\begin{enumerate}
\item each $C_\a$ is a closed cofinal subset of $\a$;
\item for each limit point $\b$ of $C_\a$, $C_\b=C_\a\cap \b$;
\item the order type of each $C_\a$ is at most $\nu$.
\end{enumerate}
Condition 3 above implies that in $V$ there is no thread $C\subseteq \kappa^+$ of $\vec{C}$. Thus, $\vec{C}$ is a $\square(\kappa^+)$ sequence in $V$. 
\end{proof}

\begin{corollary} Assume $\sf{PFA}$ and suppose $\sf{Covering\ with\ Chang\ Models}$ holds. Then there is an inner model with a subcompact cardinal.
\end{corollary}

Much has been done towards establishing covering results using canonical inner models resembling $\sf{L}$. For example, the reader can consult \cite{jensen2009stacking}, \cite{MitSchSte}, \cite{StrengthuB}, \cite{sargsyan2013covering}, \cite{SchWoo}, \cite{CMIP} and \cite{UltimateL}. All of these papers except the last one deal with the \textit{short-extender-region}. The last one has conjectures in the region of supercompact cardinals. The introduction of \cite{sargsyan2013covering} contains a lengthier discussion of \rcon{main conjecture} and its role in inner model theory.

\rcon{main conjecture} connects determinacy theories with other natural frameworks. Perhaps our approach to the program of finding canonical covering structures is somewhat unconventional as the covering structures that we propose are not constructed by conventional backgrounded constructions such as $L[\vec{E}]$ constructions or $K^c$ constructions. The reader can learn more about the conventional approach from \cite{jensen2009stacking}, which is based on the existence of $K^c$. 

Our main motivation comes from the desire to isolate large models of determinacy in natural extensions of $\sf{ZFC}$. We strongly believe that the information coded into canonical sets of reals describes the \textit{core} of the universe. Thus, our approach is closer to the $\sf{Ultimate\ L}$ framework proposed by Woodin than the conventional approach. Our ideas are also heavily influenced by Steel's inter-translatability ideas that appear in \cite{GodelSteel}. 

Besides stating \rcon{main conjecture}, our goal here is to verify that it holds in \textit{hod mice}, which are canonical structures that appear in the analysis of $\H$ of models of determinacy. The reader can learn more about them by consulting \cite{hod_mice},  \cite{BSL}, \cite{hod_mice_LSA} and \cite{normalization_comparison}.  The following is the main theorem of the paper. We heavily rely on the terminology of  \cite{normalization_comparison}.

Suppose $M$ is a transitive model of some fragment of set theory and $\k$ is a cardinal of $M$. Then $\Gamma^\k$ is the set of all $\k$-universally Baire set of reals of $M$, and if $g$ is $M$-generic then $\Gamma_g^\k$ is the set of $\k$-universally Baire sets of $M[g]$. 

Suppose $\V\models \sf{ZFC}$ is a hod premouse and  $\k$ is a limit of Woodin cardinals of $\V$. Let $g\subseteq Coll(\omega, <\k)$ be $\V$-generic. Set $g_\a=g\cap Coll(\omega, <\a)$, $\mathbb{R}^*_g = \bigcup_{\alpha<\k} \mathbb{R}^{V[g_\a]}$ and
\begin{center}
$\Gamma^*_g=\{A^g\cap \bR^*: \exists \a<\k (A\in \Gamma^\k_{g_\a})\}$.
\end{center}
The following is our main theorem.

\begin{theorem}\label{main theorem} Suppose $(\V, \Omega)$ is an excellent pair (see \rdef{excellent pair}) and $\V\models \sf{ZFC}$. Suppose further that $\kappa$ is a limit of Woodin cardinals of $\V$ such that either $\k$ is regular or its cofinality is not a measurable cardinal. Then the following stronger version of \rcon{main conjecture} is true in $\V$ at $\k$. 

There is a transitive lbr hod premouse $\M$ such that whenever $g\subseteq Coll(\omega, <\k)$ is generic. 
\begin{enumerate}
\item $Ord\cap \M=\k^+$,
\item $\M$ has a largest cardinal $\nu$,
\item for any generic $g\subseteq Coll(\omega, <\k)$, 
\begin{center}
$L(\M, \bigcup_{\a<\nu}(\M|\a)^\omega, \Gamma^*_g, \bR^*_g)\models \sf{AD}$.
\end{center}
\item If in addition $\k$ is not subcompact then $\M\models \square_\nu$.
\end{enumerate}
\end{theorem}

In fact we will give a precise definition of $\M$ as above (see \rthm{determinacy in chang model}). 

The reader may wonder if there are hod mice in which there are Woodin cardinals and strong cardinals below them. Recently, the author constructed, assuming large cardinals, hod mice that do have Woodin cardinals that are limits of Woodin cardinals. At the time of writing the current paper, this is the best known existence result on hod mice. In particular, since we present \rthm{main theorem} as an evidence for \rcon{main conjecture}, perhaps it is wiser to strengthen the minimality assumption used in \rcon{main conjecture} to the non-existence of inner models with a Woodin cardinal that is a limit of Woodin cardinals. We confess that at the time of writing the current paper, the universe beyond a Woodin cardinal that is a limit of Woodin cardinals is a somewhat  darker region. 

Another mystery may be the requirement in \rcon{main conjecture} that $\k$ be a limit of strong cardinals. After all no such requirement was made in \rthm{main theorem}. One reason for this requirement is rooted in technical arguments that have been used in our attempts to prove instances of \rcon{main conjecture}. For example, the reader may consult the introduction of \cite{StrengthuB} where proving \rcon{main conjecture} is presented as the main goal of the \textit{core model induction}. 

Another reason is that recent unpublished results of the author and Paul Larson seem to suggest that the failure of square at $\omega_3$ along with $2^{\omega}=\omega_2$ and $2^{\omega_2}=\omega_3$ is weaker than a Woodin cardinal that is a limit of Woodin cardinals. This result seems to suggest that \rcon{main conjecture} for successor cardinals or cardinals $\k$ for which $V_\k$ is not rich maybe false in that the model $M$ may not have a largest cardinal. This is exactly what seems to happen in the model mentioned above in which $\neg \square_{\omega_3}$ holds. 

Another issue with stating such covering conjectures for non-measurable regular cardinals is that for a regular cardinal $\kappa$ one can collapse $\kappa^+$ to $\kappa$ without adding bounded subsets of $\k$. This suggests that one cannot hope to obtain $M$ as in \rcon{main conjecture}, as such an $M$ is most likely absolute between $V$ and $V[g]$ where $g$ collapses $\kappa^+$ to $\kappa$ (the problematic condition seems to be the requirement that $M$ has height $\k^+$). 

Inner model theory is perhaps closer to physics than other branches of mathematics. The fundamental goal of the subject seems to be finding truths about the universe of sets rather than solving particular precisely stated problems about concrete mathematical structures or discovering truths about concrete mathematical structures such as $\mathbb{N}$ or $L$. Just like many theories describing the nature of the physical universe, inner model theory has gone through many stages of revising truths believed by the community. We certainly do not think that \rcon{main conjecture} will be the last word one will say about the short extender region of the universe of sets. There seem to be a lot of mysteries that yet to be understood, and so the intention of \rcon{main conjecture} is to suggest a concrete direction for practitioners of inner model theory. It is very likely that as we know more, just like with many theories put forward by physicists, \rcon{main conjecture} will be revisited and modified. 

\textbf{Acknowledgments.} The author wishes to thank the referee for several useful comments, especially for not trusting an earlier version of this paper. The author's work was partially supported by the NSF Career Award DMS-1352034.

\section{Terminological Interlude}\label{terminological interlude}

In this paper, our goal is to give a proof of \rthm{main theorem} without getting into nuts and bolts of inner model theory. We will thus isolate some properties of iteration strategies that Steel's \textit{mouse pairs} have and work with strategies that have these properties. We, however, will not explain how one obtains such strategies as this is simply beyond the scope of this paper (appropriate references will be given). 

In this paper, we will use some aspects of Jensen's terminology, though not so much his indexing scheme. In this paper, the indexing issues are irrelevant. However, we warn the readers that in this paper, following Jensen,  we simply use ``iterations" as opposed to ``iteration trees". \\

\textbf{A notational digression.} Suppose $\M$ is an lbr hod premouse and $\Lambda$ is an iteration strategy for $\M$. 
\begin{enumerate}
\item For any transitive set $o(M)=M\cap Ord$.
\item Given $\xi\leq Ord\cap \M$, we let $\Lambda_{\M|\xi}$ be the strategy of $\M|\xi$ induced by $\Lambda$\footnote{One can think of $\Lambda_{\M|\xi}$ as $id$-pullback of $\Lambda$.}. 
\item If $\T$ is an iteration of $\M$ according to $\Lambda$ with last model $\Q$, we let $\Lambda_{\T, \Q}$ be the strategy of $\Q$ induced by $\Lambda$. More precisely, $\Lambda_{\T, \Q}(\U)=\Lambda(\T^\frown \U)$. 
\item Suppose $\T$ is an iteration of some $\M|\k$ where $\k$ is an inaccessible cardinal of $\M$ with the property that $\rho_\omega(\M)>\k$. We then let $\T^\M$ be the $id$-copy of $\T$ on $\M$. Because of our choice of $\k$, $\T^\M$ and $\T$ have the same extenders and the same tree structure. 
\item Similarly, if $\T$ is an iteration of $\M$ and is based on $\M|\k$ where $\k$ is an inaccessible cardinal of $\M$ with the property that $\rho_\omega(\M)>\k$, we let $\T\rest \M|\k$ be the iteration of $\M|\k$ obtained by (essentially) copying $\T$ on $\M|\k$\footnote{The fact that $\k$ is inaccessible or that $\rho_\omega(\M)>\k$ are largely irrelevant for these two notions, provided $\M$ has an iteration strategy with strong hull condensation.  However defining the meaning of $\T\rest \M|\k$ and $\T^\M$ when $\kappa$ is arbitrary is not as straightforward. We do not need them in this paper, so we will not be dealing with these notions.}. 
\item Given an interval $(\a, \b)$ with $\b\leq o(\M)$, we say that iteration $\T$ is based on $\M|(\a, \b)$ or just $(\a, \b)$ if all extenders used in $\T$ have critical points $>\a$ and for each $\gg<lh(\T)$, if $\pi^\T_{0, \gg}$ exists then $lh(E_\gg)<\pi^\T_{0, \gg}(\b)$. We say $\T$ is based on $\M|\b$ if it is based on $(0, \b)$.
\item We say $\T$ is above $\a$ if $\T$ is based on $(\a, o(\M))$. 
\item If $\T$ is an iteration of $\M$ and $\a<lh(\T)$ then $\T_{\leq \a}=\T\rest \a+1$ and $\T_{\geq \a}$ is $\T$ after stage $\a$. 
\item Suppose $\Q$ is a $\Lambda$-iterate of $\M$ and $\R$ is a $\Lambda_\Q$-iterate of $\Q$. Suppose $\Lambda$ has \textit{full normalization} (see clause 1 of \rdef{excellent pair}). Then we let $\T^\Lambda_{\Q, \R}$ be the unique normal $\Lambda_\Q$-iteration of $\Q$ with last model $\R$. If $\pi^{\T^\Lambda_{\Q, \R}}$ exists then we let $\pi^\Lambda_{\Q, \R}$ be this embedding. When $\Lambda$ is clear from context, we will drop it from our notation.
\item We say that $(\N, \Phi)$ is an iterate of $(\M, \Lambda)$ if $\N$ is a $\Lambda$-iterate of $\M$ via $\T$ and $\Phi=\Lambda_{\T, \N}$. 
\item We say $\T$ is strongly non-dropping if for each $\a<lh(\T)$, $\nu_\a^\T$, the supremum of the generators of $E_\a^\T$, is an inaccessible cardinal in $\M_\a^\T$. 
\item We let $S^\M$ is the strategy predicate of $\M$.
\item In this paper, to keep the matters simple, when we say that some cardinal has large cardinal properties in a fine structural model we always mean that the large cardinal properties in question are witnessed by the extenders on the sequence. Schlutzenberg has done substantial work showing that many such large cardinals are indeed witnessed by the extenders on the sequence of the fine structural model (see for example \cite{extmax}). 
\item  We say $w=(\eta^w, \d^w)$ is a \textit{window} of $\M$  if $\M\models ``$there are no Woodin cardinals in the interval $(\eta^w, \d^w)$ and $\d^w$ is a Woodin cardinal". 
\item We say $w$ is a \textit{maximal window} if in addition $\eta^w$ is the least $\M$-inaccessible $>\sup\{\xi: \M|\d^w\models ``\xi$ is a Woodin cardinal"$\}$.
\item Given two windows $w$ and $v$, we write $w<_W v$ if $\d^w<\d^v$ (notice that this implies that $\d^w<\eta^v$). 
\item We say that an iteration $\T$ of $\M$ is based on $w$ if $\T$ is based on $(\eta^w, \d^w)$.
\item We let ${\sf{EA}}^\M_w$ be the extender algebra of $\M$ associated with $\d^w$ that only uses extenders $E$ such that $\cp(E)>\eta^w$ and $\nu_E$, the sup of the generators of $E$, is an inaccessible cardinal of $\M$\footnote{We will not specify which generator version of the extender algebra we use, though we will never use $\d^w$-generator version.}.
\end{enumerate}

\section{Excellent hod pairs}

In this section, we isolate some abstract properties of Steel's notion of \textit{mouse pairs} (see \cite[Chapter 0.1]{MousePairs}). These abstract properties are the properties that we will need in our proof of \rthm{main theorem}.

\begin{definition}\label{excellent pair} Suppose $\Sigma$ is a strategy for an lbr hod mouse $\P$. We say that $\Sigma$ is \textit{almost excellent} if
\begin{enumerate}
\item $\Sigma$ has \textbf{full normalization}, i.e., whenever $\T$ is an iteration of $\P$ according to $\Sigma$ with last model $\Q$, there is a normal iteration $\U$ of $\P$ according to $\Sigma$ with last model $\Q$ and such that 
\begin{enumerate}
\item $\pi^\T$ exists if and only if $\pi^\U$ exists,
\item if $\pi^\T$ exists then $\pi^\T=\pi^\U$, and
\item if $\T$ is strongly non-dropping then $\U$ is strongly non-dropping.
\end{enumerate} 
\item $\Sigma$ is \textbf{positional}, i.e., if $\Q$ is a $\Sigma$-iterate of $\P$ via both $\T$ and $\U$ then $\Sigma_{\T, \Q}=\Sigma_{\U, \Q}$\footnote{We will thus drop $\T$ from $\Sigma_{\T, \Q}$.},
\item $\Sigma$ is \textbf{stable}, i.e., if $\Q$ and $\R$ are $\Sigma$-iterates of $\P$ then $\Q\insegeq \R$ implies that $\Q=\R$,
\item $\Sigma$ is \textbf{segmentally normal}, i.e., whenever $\eta$ is an inaccessible cardinal of $\P$ such that $\rho_{\omega}(\P)>\eta$, $\T$ is a strongly non-dropping $\Sigma$-iteration of $\P$ that is above $\eta$, $\Q$ is the last model of $\T$,  $\U$ is a $\Sigma_\Q$-iteration of $\Q$ that is based on $\Q|\eta$ such that $\pi^\U$ exists and $\R$ is the last model of $\U$ then 
\begin{enumerate}
\item $\Sigma_{\P|\eta}=(\Sigma_{\Q})_{\P|\eta}$  and
\item letting $\W=(\U\rest (\Q|\eta))^\P$\footnote{See clause 4 and 5 of \rsec{terminological interlude}.} and $\S$ be the last model of $\W$, $\R$ is a normal $\Sigma_\S$-iterate of $\S$ that is strongly non-dropping and is above $\pi_{\P,\S}(\eta)$.
\end{enumerate}
%
\item $\Sigma$ is \textbf{directed}, i.e., if $\Q$ and $\R$ are $\Sigma$-iterates of $\P$ then there is $\S$ that is a $\Sigma_\Q$-iterate of $\Q$ and a $\Sigma_\R$-iterate of $\R$, and $\card{\S}=\max(\card{\Q}, \card{\R})$.
\item $(\P, \Sigma)$ satisfies \textbf{generic interpretability}, i.e., for any $\k$ that is a limit of Woodin cardinals of $\P$ and for any $\eta<\k$ and any $\P$-generic $g\subseteq Coll(\omega, \eta)$, setting $\Lambda=\Sigma\rest (HC^{\P[g]})$,
\begin{enumerate}
\item $\Lambda \in \P[g]$,
\item $\Lambda \rest (HC^{\P[g]})$ is $\k$-universally Baire in $\P[g]$ as witnessed by trees $T, S$ that are definable over $\P[g]$ from the pair $(g, S^\P_{\P|\eta})$,
\item if $(T, S)\in \P[g]$ is any pair of trees witnessing clause 2 above, for any $<\k$-generic $h$ over $\P[g]$, $(p[T])^{\P[g][h]}=\Sigma\rest (HC^{\P[g][h]})$,
\item in $\P[g]$, $\Lambda$ is the unique strategy $\Phi$ of $\P|\eta$ that extends $S^\P_{\P|\eta}\rest (HC^{\P[g]})$ and satisfies clauses (b).
\end{enumerate}
\item $\Sigma$ is \textbf{pullback consistent}, i.e., if $\Q$ is a $\Sigma$-iterate of $\P$ such that $\pi_{\P, \Q}$ is defined and $\eta<Ord\cap \P$ then $\Sigma_{\P|\eta}$ is the $\pi$-pullback of $\Sigma_{\Q|\pi_{\P, \Q}(\eta)}$. 
\end{enumerate} 
We say that $\Sigma$ is \textbf{excellent} if whenever $\M\insegeq \P$ is such that $o(\M)$ is an inaccessible cardinal of $\P$ and $\rho_\omega(\P)>o(\M)$, $\Sigma_\M$, the $id$-pullback of $\Sigma$, is almost excellent. 

In this paper, we say $(\P, \Sigma)$ is an \textbf{excellent pair} if $\P$ is a countable lbr hod premouse and $\Sigma$ is an excellent $(\omega_1, \omega_1+1)$-iteration strategy for $\P$. 
\end{definition}

The reader might be wondering if $\Lambda$ of clause 6 above is directed in $\P[g]$. This is an easy consequence of genericity iterations. First, we will let $\Lambda^g$ be the unique $\Lambda$ as in clause 6, and $\Lambda^{g*h}$ be its canonical extension in $\P[g][h]$.

\begin{lemma} Suppose $(\P, \Sigma)$ is an excellent pair, $\k$ is a limit of Woodin cardinals of $\P$ and $\eta<\k$. Let $g\subseteq Coll(\omega, \eta)$ be $\P$-generic and $h\subseteq Coll(\omega, <\k)$ be generic over $\P[g]$. Then $\Lambda^{g*h}\rest HC^{\P(\bR^*_{g*h})}$ is directed.  
\end{lemma}
\begin{proof} Let $\Q$ and $\R$ be $\Lambda^{g*h}$-iterates of $\P|\eta$ in $\P(\bR^*_{g*h})$. Let  $\xi\in (\eta, \k)$ be such that $\Q, \R\in \P|\xi[g*h_\xi]$, and let $\d\in (\xi, \k)$ be a Woodin cardinal of $\P[g*h_\xi]$ and $\S$ be a common iterate of $\Q$ and $\R$ via $\Sigma_{\Q}$ and $\Sigma_{\R}$. Let $\P_0$ be a non-dropping iteration of $\P$ based on $(\xi, \d)$ such that $\S$ is generic over $\P_0$. Then  $\P_0[g*h_\xi][\S]$ thinks that $\Q$ and $\R$ have a common $\Lambda^{g*h_\xi}_\Q$ and $\Lambda^{g*h_\xi}_\R$ iterate. This fact pulls back to $\P[g*h_\xi]$, and as $\P|\d$ is countable in $\P[g*h]$, we can find $\S \in HC^{\P(\bR^*_{g*h})}$ that is a common $\Lambda^{g*h}_\Q$ and $\Lambda^{g*h}_\R$ iterate of $\Q$ and $\R$ . 
\end{proof}

\begin{remark} 
\begin{enumerate}
\item The assumptions made in \rdef{excellent pair}, while legitimate, are not part of the defining conditions of a hod premouse; they are only  consequences of Steel's notion of \textbf{mouse pair} which, for example, appears in \cite[Chapter 0.1]{MousePairs}. We advice the reader unfamiliar with the nuts and bolts of normalization to glance over \cite[Chapter 0.1]{MousePairs}, and consider the references mentioned in that chapter. However, exactly how various normalization procedures work are irrelevant for the current paper. 

\item The following are the relevant definitions and theorems of \cite{MousePairs} that the reader is encouraged to consult: Definition 0.1, Theorem 0.3 and the entire Chapter 1. These definitions and results have their origins in \cite{normalization_comparison}. The reader can easily locate them in \cite{normalization_comparison} by following the references listed in \cite{MousePairs}.  

\item Notice that clause 1a, 1b and 2 of \rdef{excellent pair} are simply Theorem 1.1 and Corollary 1.2 of \cite{MousePairs}. Clause 3 is a consequence of the Dodd-Jensen property (see Theorem 0.3 of \cite{MousePairs}). Clause 1c is a consequence of embedding normalization (see the next comment).

\item Clause 4 is an immediate consequence of embedding normalization as spelled out in \cite[Chapter 3]{normalization_comparison}. For example, consider \cite[Claim 3.2]{normalization_comparison} and \cite[Proposition 3.18]{normalization_comparison}.

\item Clause 5 is a consequence of the comparison theorem (see \cite[Theorem 1.3]{{normalization_comparison}}).
\item Clause 6 is verified in \cite[Chapter 8.1, Theorem 8.1]{normalization_comparison}.
\item Clause 7 is a standard consequence of \textit{hull condensation} and can be verified by consulting \cite[Lemma 4.9]{normalization_comparison}.

\item The requirement that $\Sigma_\M$ is almost excellent is a consequence of the fact that initial segments of mouse pairs are mouse pairs (see \cite{normalization_comparison}). However, notice that this requirement is made only for those $\M$ for which $o(\M)$ is an inaccessible cardinal of $\P$ and $\rho_\omega(\P)>o(\M)$. Given such an $\M$ and an iteration $\T$ of it, there is essentially no difference between $\T^\P$ and $\T$. Thus, for such $\M$, $\Sigma_\M$ straightforwardly inherits all of the condensation properties that $\Sigma$ might have.

\item It is straightforward to verify that if $\Sigma$ is excellent and $\Q$ is a $\Sigma$-iterate of $\P$ then $\Sigma_\Q$ is also excellent.
\end{enumerate}
\end{remark}

Thus, the following theorem, due to Steel, summarizes the above remark.

\begin{theorem} Assume $AD^+$ and suppose $(\P, \Sigma)$ is a mouse pair such that $\P$ is an lbr hod premouse. Then $(\P, \Sigma)$ is an excellent pair.
\end{theorem}

If $\Sigma$ is an excellent strategy and $\Q$ is a $\Sigma$-iterate of $\P$ via $\T$ then $\Sigma_{\T, \Q}$ is independent of $\T$, and so we will drop $\T$ from our notation (see \cite[Corollary 1.2]{MousePairs}). 

The following is a simple consequence of excellence that we will need in this paper. 

\begin{lemma}\label{initial segment normalization} Suppose $(\P, \Sigma)$ is an excellent pair and 
\begin{enumerate}
\item $\lambda$ is an inaccessible cardinal of $\P$ such that $\rho_\omega(\P)>\l$, and
\item $\Q_0$ and $\Q_1$ are $\Sigma$-iterates of $\P$ via normal iterations $\U_0$ and $\U_1$ respectively such that both $\U_0$ and $\U_1$ are based on $\P|\l$ and do not have drops on their main branch.
\end{enumerate}
Let for $k\in 2$, $\l_k=\pi^{\U_k}(\l)$. Suppose next that $(\S, \Phi)$ is an iterate of both $(\Q_0|\l_0, \Sigma_{\Q_0|\l_0})$ and $(\Q_1|\l_1, \Sigma_{\Q_1|\l_1})$. Let\footnote{See the above notational digression.}
\begin{center}
$\U_2=(\T_{\Q_0|\l_0, \S})^{\Q_0}$ and $\U_3=(\T_{\Q_1|\l_1, \S})^{\Q_3}$.
\end{center}
 Finally, let $\Q_2$ be the last model of $\U_2$ and $\Q_3$ be the last model of $\U_3$. Then $\Q_2=\Q_3$ (and clearly, $\Sigma_{\Q_2}=\Sigma_{\Q_3}$
).
\end{lemma}
%
\begin{proof} First note that both $\pi^{\U_2}$ and $\pi^{\U_3}$ are defined and setting $\nu=\S\cap Ord$,
\begin{center}
$\pi^{\U_2}(\l_0)=\nu=\pi^{\U_3}(\l_1)$.
\end{center}
The above equality is a consequence of the full normalization of $\Sigma_{\P|\l}$. Indeed, full normalization of $\Sigma_{\P|\l}$ implies that
\begin{center}
$\pi^{\U_2\rest \Q_0|\l_0}\circ \pi^{\U_0} =\pi^{\U_3\rest \Q_1|\l_1}\circ \pi^{\U_1}$.
\end{center}
Because $\l$ is inaccessible, we have that both $\pi^{\U_2}\circ \pi^{\U_0}$ and $\pi^{\U_3}\circ \pi^{\U_1}$ are continuous at $\l$. 

Let now $\T$ be the normal iteration of $\P|\l$ according to $\Sigma_{\P|\l}$ with last model $\S$. Thus, $\T$ is the full normalization of both $(\U_0\rest \P|\l)^\frown (\U_2\rest \Q_0|\l_0)$ and $(\U_1\rest \P|\l)^\frown (\U_3\rest \Q_1|\l_1)$, and it is obtained via the least extender comparison process between $\P|\l$ and $\S$. Because $\Sigma$ is excellent, the last branch of $\T$ doesn't drop. An argument like the one given above implies that $\pi^{\T^\P}(\l)=\nu$. Let $\R$ be the last model of $\T^\P$. We have that $\S\insegeq \R$. 

We now claim that $\R=\Q_2$. Suppose not. Clearly, $\R|\nu=\Q_2|\nu=\S$. Let $\W$ be the normal iteration of $\P$ according to $\Sigma$ with last model $\Q_2$. Thus, $\W$ is the normalization of $\U_0^\frown \U_2$, and we must have that $\pi^\W=\pi^{\U_2}\circ \pi^{\U_0}$. Because we are assuming $\R\not =\Q_2$, it follows that $\T\inseg \W$ (recall that $\W$ is just the iteration produced by the least extender comparison between $\P$ and $\Q_2$ and $\P$-to-$\R$ iteration is part of it). Let $\a$ be such that $\M_\a^\W=\R$. We then must have that $lh(E_\a^\W)>\nu$\footnote{Equality cannot hold because $\nu$ is inaccessible in $\R$ and $<$ cannot hold because otherwise we will have $\R|\nu\not =\S$.}.

Because $\U_0^\frown \U_2$ is based on $\P|\l$ we have that the generators of $\U_0^\frown \U_2$ are contained in $\nu$. 
Thus, the generators of $\W$ are contained in $\nu$. Notice that if $\cp(E_\a^\W)<\nu$ then we cannot have that $\nu\in rng(\pi^\W)$ which contradicts full normalization as $\nu\in rng(\pi^{\U_2}\circ \pi^{\U_0})$\footnote{In this case if $\nu\in rng(\pi^\W)$ then we must have that $\a$ is on the main branch of $\W$ and the extender used at $\a$ on the main branch has a critical point in the interval $(\nu, lh(E_\a^\W))$. This easily implies that $\W$ must have generators above $\nu$.}. If $\cp(E_\a^\W)>\nu$ then $\W$ has generators above $\nu$ implying that $\U_0^\frown \U_2$ must also have such a generator.

By a symmetric argument, we also have that $\R=\Q_3$. Hence, $\Q_2=\Q_3$.
\end{proof}

%
%
%

\section{Internal direct limit constructions}\label{internal dir lim cons: sec}

Suppose $(\V, \Omega)$ is an excellent pair and $\k$ is a limit of Woodin cardinals of $\V$ such that if $\cf^\V(\k)<\k$ then $\cf^\V(\k)$ is not a measurable cardinal of $\V$. We will drop $\Omega$ from now on. In fact, its only use is to guarantee that the internal strategy of $\V$ has the properties stated in \rdef{excellent pair}. 

 Set $\P=\V|(\k^+)^\V$. We let $\Sigma$ be the $(\k, \k+1)$-fragment\footnote{I.e., $\k$-rounds each of which can be $\k+1$ length. Player I can start a new round only if the previous round has length strictly shorter than $\k$.} of $S^\V_\P$ that acts on iterations that are based on $\P|\k$.  Given any $\V$-generic $h$, we let $\Sigma^h$ be the extension of $\Sigma$ in $\V[h]$ (see \rdef{excellent pair})\footnote{We will only use this notation for posets of size strictly smaller than $\k$ or for $Coll(\omega, <\k)$.}. Suppose $g\subseteq Coll(\omega, <\k)$ is $\V$-generic. The following is a notation that we will use in this paper.
\begin{enumerate}
 \item For $\a<\k$, let $g_\a=g\cap Coll(\omega, <\a)$. 
 \item Let $\mathcal{I}^{g_\a}(\P)$ be the set of $\Sigma^{g_\a}$-iterates of $\P$ that are obtained via iterations $\T\in V[g_\a]$ such that $lh(\T)\leq \k+1$, $\T$ is based on $\P|\k$, $\pi^\T$ is defined and $\pi^\T(\k)=\k$. 
 \item Set $\mathcal{I}^g(\P)=\cup_{\a<\k}\mathcal{I}^{g_\a}(\P)$.
\item Given $\Q\in \mathcal{I}^{g_\a}(\P)$ and $\b\in [\a, \k)$, we let $\mathcal{F}^{g_\b}_\Q$ be the set of $\Sigma^{g_\b}_\Q$-iterates $\R$ of $\Q$ such that $lh(\T_{\Q, \R})<\k$ and $\pi_{\Q, \R}$ is defined. 
\item Set $\mathcal{F}^{g,\a}_\Q=\cup_{\b\in (\a, \k)}\mathcal{F}^{g_\b}_\Q$. Clearly, $\mathcal{F}^{g, \a}_\Q=\mathcal{F}^{g, \b}_\Q$, and so we drop $\a$ from our notation and just write $\mathcal{F}^{g}_\Q$. 
\end{enumerate}

Below we introduce window-based iterations and genericity-iterations.

\begin{definition}\label{window-based iterations}
Suppose $\R\in \mathcal{I}^g(\P)$. We say $\Q$ is a \textit{window-based iterate} of $\R$ if there is $\iota<\k$ such that $\R\in V[g_\iota]$, an $<_W$-increasing sequence of windows $(w_\a=(\eta_\a, \d_\a): \a<\cf(\k))$ of $\R$ and a sequence $(\Q_\a: \a<\cf(\k))\subseteq \mathcal{F}^{g_\iota}_\R$ (in $V[g_\iota]$) such that
\begin{enumerate}
\item for all $\a<\cf(\k)$, $\d_\a<\k$,
\item $\sup_{\a<\cf(\k)}\d_\a=\k$,
\item $\Q_0\in \mathcal{F}^{g_\iota}(\R)$ and $\T_{\R, \Q_0}$ is based on $\R|\eta_0$,
\item $\Q_{\a+1}\in \mathcal{F}^{g_\iota}_{\Q_\a}$,
\item $\Q_{\a+1}$ is obtained from $\Q_\a$ via an iteration according to $\Sigma^{g_\iota}_{\Q_\a}$ that is based on $\pi_{\R, \Q_\a}(w_\a)$,
\item for limit ordinals $\l$, $\Q_\l$ is the direct limit of $(\Q_\a, \pi_{\Q_\a, \Q_\b}: \a<\b<\l)$,
\item $\Q$ is the direct limit of $(\Q_\a, \pi_{\Q_\a, \Q_\b}: \a<\b<\cf(\k))$.
\end{enumerate}
%
We set $\eta_{\a, \b}=\pi_{\R, \Q_\a}(\eta_\b)$, $\eta^\a=\eta_{\a, \a}$, $\d_{\a, \b}=\pi_{\R, \Q_\a}(\d_\b)$ and $\d^\a=\d_{\a, \a}$. 
\end{definition}

It should be clear to the reader that there are unique sequences  $(w_\a=(\eta_\a, \d_\a): \a<\cf(\k))$ and $(\Q_\a: \a<\cf(\k))$ witnessing that $\Q$ is a window-based iterate of $\R$.

\begin{definition}\label{genericity iterations}
Suppose $\R\in \mathcal{I}^g(\P)$. We say that $\Q$ is a \textit{genericity iterate} of $\R$ if it is a window-based iterate of $\R$ as witnessed by $(w_\a: \a<\cf(\k))$ and $(\Q_\a: \a<\cf(\k))$ such that
\begin{enumerate}
\item if $x\in \bR_g$ then for some $\a<\k$, $x$ is generic for the extender algebra ${\sf{EA}}^\Q_{\pi_{\R, \Q}(w_\a)}$,
\item for each $\a<\cf(\k)$, $w_\a\in rng(\pi_{\P, \R})$,
\item for each $\a$, $\T_{\Q_\a, \Q_{\a+1}}$ is a strongly non-dropping iteration of $\Q_\a$.
\end{enumerate} 
\end{definition}

Clause 1 of \rdef{genericity iterations} may seem odd to the reader as $x\in \bR_g$ while $\Q$ is a $\Sigma^{g_\iota}_\R$-iterate of $\R$ for some $\iota<\k$. Fix an inaccessible cardinal $\nu\in (\iota, \k)$ such that $x\in \P[g_\nu]$. In $\V[g_\iota]$, given a window $w$ of $\R$, we can iterate $\R$ in the window $w$ to make $a=_{def}\P|\nu^+[g_\iota]$ generic. If $\S$ is the result of this iteration then $x$ is generic over $\S$\footnote{Notice that $\T_{\P, \S}\in \P[g_\iota]$. Thus, $H_{\nu^+}^{\S[a]}=\P|\nu^+[g_\iota]$. Since $x$ is generic over $\P|\nu^+[g_\iota]$, $x$ is generic over $\S[a]$, and hence, over $\S$.}. If $\Q$ now is an iterate of $\R$ such that for each $\nu<\k$, $\P|\nu^+$ is generic over some window of $\Q$ then clause 1 of \rdef{genericity iterations} will be satisfied. Notice that because we are using $\sf{EA}$ as our extender algebra, the genericity iterations described above are all strongly non-dropping. The following, then, follows from the above discussion.

\begin{proposition}\label{constructing genericity iterations} Suppose $\R$ is a genericity iterate of $\P$ and $\S\in \mathcal{F}^g_\R$. Then $\T_{\R, \S}$ can be normally continued to a genericity iterate of $\R$. More precisely, there is $\Q$ that is a genericity iterate of $\S$ such that $(\T_{\R, \S})^\frown \T_{\S, \Q}$ is normal and is a genericity iterate of $\R$.  
\end{proposition}

The following is a consequence of excellence.

\begin{proposition}\label{comute genericity} Suppose $\R$ is a genericity iterate of $\P$, $\Q$ is a genericity iterate of $\R$ and $\S$ is a genericity iterate of $\Q$. Then $\S$ is a genericity iterate of $\R$. 
\end{proposition}
\begin{proof} Let $(w_\a: a<\cf(\k))$ be the windows used in $\T_{\R, \Q}$ and $(u_\a: \a<\cf(\k))$ be the windows used in $\T_{\Q, \S}$. Notice that the full normalization of $(\T_{\R, \Q})^\frown \T_{\Q, \S}$ is strongly non-dropping. For each $\a<\cf(\k)$, set $u_\a'=\pi_{\R, \Q}^{-1}(u_\a)$. Notice that for each $\a<\cf(\k)$, $u_\a'\in rng(\pi_{\P, \R})$. Let $(v_\a:\a<\cf(\k))$ be the enumeration of $\{w_\a:\a<\cf(\k)\}\cup \{u_\a': \a<\cf(\k)\}$ in increasing order. It then follows from excellence that $\S$ is a window iterate of $\R$ as witnessed by $(v_\a:\a<\cf(\k))$. It remains to see that for each $x$ there is $\a<\cf(\k)$ such that $x$ is generic for ${\sf{EA}}^\S_{\pi_{\R, \S}(v_\a)}$. Let $\a$ be such that  $x$ is generic for ${\sf{EA}}^\S_{\pi_{\Q, \S}(u_\a)}$ and let $\b$ be such that $u_\a'=v_\b$. Then $x$ is generic for ${\sf{EA}}^\S_{\pi_{\R, \S}(v_\b)}$.
\end{proof}

The following is a straightforward consequence of excellence and \rlem{initial segment normalization}.
\begin{lemma}\label{directedness} Suppose $\Q$ is a window based iterate of $\P$ and $\R_0, \R_1\in \mathcal{F}^g_\Q$. There is then an $\R\in \mathcal{F}^g_\Q$ such that $(\R, \Sigma^g_\R)$ is a common iterate of $(\R_0, \Sigma^g_{\R_0})$ and $(\R_1, \Sigma^g_{\R_1})$.
\end{lemma}
\begin{proof}
Let $\T=\T_{\P, \Q}$, $\U_0=\T_{\Q, \R_0}$ and $\U_1=\T_{\Q, \R_1}$. Notice that because for $i\in2$, $lh(\U_i)<\k$, we must have that there is some $\eta'<\k$ such that $\eta'$ is an inaccessible cardinal of $\Q$ and for $i\in 2$, $\U_i$ is a normal iteration of $\Q$ based on $\Q|\eta'$. Let $\eta$ be the least such $\eta'$ and let $\xi<lh(\T)$ be such that $\M_\xi^\T|\eta=\Q|\eta$. Let $\l_0$ be the least inaccessible of $\P$ such that $\T\rest \xi$ is based on $\P|\l_0$. Let $(w_\a=(\eta_\a, \d_\a): \a<\cf(\k))$ be the windows used in $\T$. Let $\l=\eta_\a$ where $\a$ is the least such that $\l_0<\eta_\a$. Because $\T$ is a window-based iteration, it follows that if $\zeta<\k$ is the least such that $\M_\zeta^\T|\pi^\T(\l)=\Q|\pi^\T(\l)$ then (notice that $\pi^\T(\l)=\eta^\a$) \\\\
(a)  $\T_{\geq \zeta}$ is a normal window-based iteration of $\M_\zeta^\T$ that is above $\pi^\T(\l)$, and\\
(b) $\T\rest \zeta+1$ is based on $\P|\l$. \\\\  
Set $\nu=\pi^\T(\l)$. Let now for $k\in 2$, $\S_k$ be the last model of $\U_k\rest \Q|\nu$. Using the fact that $\Sigma^g_{\P|\l}$ is directed (see \rdef{excellent pair}) we can find a common iterate $(\S, \Phi)$ of $(\S_0, \Sigma^g_{\S_0})$ and $(\S_1, \Sigma^g_{\S_1})$. Let $\R$ be the last model of $(\T_{\S_0, \S})^{\R_0}$. It follows from \rlem{initial segment normalization} that $\R$ is the last model of $(\T_{\S_1, \S})^{\R_1}$. 
\end{proof}

Suppose now that $\Q$ is a window based $\Sigma^g$-iterate of $\P$. For $\R, \S\in \mathcal{F}^g_\Q$, we let $\R\leq_\Q \S$ if $\S\in \mathcal{F}^g_\R$. It follows from \rlem{directedness} that $\leq_\Q$ is directed. Moreover, notice that if $\R, \S\in \mathcal{F}^g_\Q$ and $\W\in \mathcal{F}^g_\R\cap \mathcal{F}^g_\S$ then $\W$ is a normal $\Sigma^g_\Q$-iterate of $\Q$  and
\begin{center}
$\pi_{\R, \W}\circ \pi_{\Q, \R}=\pi_{\Q, \W}=\pi_{\S, \W}\circ \pi_{\Q, \S}$. 
\end{center}
We can then let $\M_\infty(\Q)$ be the direct limited of the directed system $(\mathcal{F}^g_\Q, \leq_\Q)$ under the iteration embedding. Given $\R\in \mathcal{F}^g_\Q$ we let $\pi^\Q_{\R, \infty}: \R\rightarrow \M_\infty(\Q)$ be the direct limit embedding. We also set $\k^\Q_\infty=\pi_{\Q, \infty}^\Q(\k)$.

Suppose $\Q$ is a genericity iterate of $\R$. There is then $h\subseteq Coll(\omega, <\k)$ with $h\in \V[g]$ such that $h$ is $\Q$-generic and $(\bR^*_h)^{\Q[h]}=\bR^*_g$. We call such an $h$ a \textit{maximal generic}. 

\begin{lemma}\label{-1} Suppose $\Q$ is a genericity iterate of $\P$. Let $h\subseteq Coll(\omega, <\k)$ be a maximal $\Q$-generic. Then $\M_\infty(\Q)=(\M_\infty(\Q))^{\Q[h]}$.
\end{lemma}
\begin{proof} This is simply because $\mathcal{F}^g_\Q=(\mathcal{F}^h_\Q)^{\Q[h]}$.
\end{proof}

Suppose $\iota<\k$ and $\Q\in \mathcal{I}^{g_\iota}(\P)$ is a window-based iterate of $\P$ as witnessed by $(w_\a=(\eta_\a, \d_\a): \a<\cf(\k))$ and $(\Q_\a: \a<\cf(\k))$. We set up the following notation.
\begin{enumerate}
\item Suppose $\R\in \mathcal{F}^g_\Q$ and $y\in \R|\k$. Let $\a_{\R, y}$ be the least $\a$ such that $\T_{\Q, \R}$ is based on $\Q|\eta^\a$ and $y\in \R|\pi_{\Q, \R}(\eta^\a)$. 
\item Given $\b\geq \a_{\R, y}$, let $\W(\R, y, \b)$ be the last model of $(\T_{\Q, \R}\rest (\Q|\eta^\a))^{\Q_\b}$. 
\end{enumerate}
The following lemma is an easy consequence of excellence (see \rdef{excellent pair}); in particular, of segmental normality. 
\begin{lemma}\label{0} Suppose $\iota<\k$ and $\Q\in \mathcal{I}^{g_\iota}(\P)$ is a genericity iterate of $\P$ as witnessed by $(w_\a=(\eta_\a, \d_\a): \a<\cf(\k))$ and $(\Q_\a: \a<\cf(\k))$. Let $\xi\in [\iota, \k)$, $\R$ and $y$ be such that $\R\in \mathcal{F}^{g_\xi}_\Q$ and $y\in \R|\k$. Let $\b \leq \gg <\cf(\k)$ be such that $\a_{\R, y}\leq \b$. Then $\W(\R, y, \gg)$ is a normal $\Sigma^{g_\xi}_{\W(\R, y, \b)}$-iterate of $\W(\R, y, \b)$ via a normal strongly non-dropping iteration  $\U$ such that $\U$ is above $\pi_{\Q_\b, \W(\R, y, \b)}(\eta^\b)$ and $\U$ is based on $\W(\R, y, \b)|\pi_{\Q_\b, \W(\R, y, \b)}(\d_{\b, \gg})$. 
\end{lemma}
\begin{proof}
First apply segmental normality to $\Q_\b|\d_{\b, \gg}$. Then apply \rlem{initial segment normalization}.
\end{proof}

The following is the main theorem of this section.

\begin{theorem}\label{1} Suppose $\Q$ is a genericity iterate of $\P$ as witnessed by $(w_\a=(\eta_\a, \d_\a): \a<\cf(\k))$ and $(\Q_\a: \a<\cf(\k))$. Then 
\begin{center}
$\M_\infty(\Q)=\M_\infty(\P)$.
\end{center}
\end{theorem}
\begin{proof} We define $j: \M_\infty(\Q)|\k_\infty^\Q\rightarrow \M_\infty(\P)|\k_\infty^\P$ as follows. Given $x\in \M_\infty(\Q)|\k_\infty^\Q$ let $\R\in \mathcal{F}^g_\Q$ be such that $\pi^{\Q}_{\R, \infty}(y)=x$ for some $y\in \R$. Notice that for every $\b\in [\a_{\R, y}, \cf(\k))$, $\W(\R, y, \b)\in \mathcal{F}^g_\P$ . Set
\begin{center}$j(x)=\pi^\P_{\W(\R, y, \a_{\R, y}), \infty}(y)$.
\end{center}
%
\textit{Claim 1.} For $\b\geq \a_{\R, y}$, $j(x)=\pi^\P_{\W(\R, y, \b), \infty}(y)$.\\\\
\begin{proof} Set $\a=\a_{\R, y}$. Because $\cp(\pi_{\W(\R, y, \a), \W(\R, y, \b)})>\pi_{\Q_\a, \W(\R, y, \a)}(\eta^\a)$, we have the following equalities.
\begin{align*}
\pi^\P_{\W(\R, y, \a), \infty}(y)&=\pi^\P_{\W(\R, y, \b), \infty}(\pi_{\W(\R, y, \a), \W(\R, y, \b)}(y))\\
&=\pi^\P_{\W(\R, y, \b), \infty}(y)
\end{align*}
\end{proof}

\textit{Claim 2.} $j$ is well-defined.\\\\
\begin{proof} Fix $\xi\in [\iota, \k)$, $\R, \R'\in \mathcal{F}^{g_\xi}_\Q$, $y\in \R$ and $y'\in \R'$ such that $\pi^\Q_{\R, \infty}(y)=\pi^\Q_{\R', \infty}(y')=x$. Let $\b>\max(\a_{\R, y}, \a_{\R', y'})$. Because of Claim 1 above, it is enough to show that 
\begin{center}
$\pi^\P_{\W(\R, y, \b), \infty}(y)=\pi^\P_{\W(\R', y', \b), \infty}(y')$.
\end{center}
Set $\W=\W(\R, y, \b)$ and $\W'=\W(\R', y', \b)$. Let $\eta=\pi_{\Q_\b, \W}(\eta^\b)$ and $\eta'=\pi_{\Q_\b, \W'}(\eta^\b)$.
Notice that $\Sigma^{g_\xi}_{\Q_\b|\eta^\b}$ is excellent.

Using the fact that $\Sigma^g_{\Q_\b|\eta^\b}$ is directed we can find a common $\Sigma^g_{\W|\eta}$ and $\Sigma^g_{\W'|\eta'}$ iterate $\M$ of $\W|\eta$ and $\W'|\eta'$. Set $\U=\T_{\W|\eta, \M}$ and $\U'=\T_{\W'|\eta', \M}$.  It follows from excellence that 
\begin{center}
$\pi_{\W|\eta, \M}\circ \pi_{\Q_\b|\eta^\b, \W|\eta}=\pi_{\W'|\eta', \M}\circ \pi_{\Q_\b|\eta^\b, \W'|\eta'}$.
\end{center}
Set 
\begin{enumerate}
\item $\Y=\T_{\Q_\b|\eta^\b, \W|\eta}$ and $\Y'=\T_{\Q_\b|\eta^\b, \W'|\eta'}$,
\item $\X=\T_{\Q_\b|\eta^\b, \M}$, $\X_0=\X^{\Q_\b}$ and $\X_1=\X^{\Q}$, 
\item $\U_0=\U^\W$ and $\U'_0=\U'^{\W'}$,
\item $\U_1=\U^\R$ and $\U_1'=\U'^{\R'}$.
\end{enumerate} 
%
%
Let $\S$ be the last model of $\X_0$.  Then $\S$ is also the last model of $\U_0$ and $\U'_0$ (see \rlem{initial segment normalization}).

Notice that $\X_1$ is the full normalization of $(\T_{\Q, \R})^\frown \U_1$ and $(\T_{\Q, \R'})^\frown \U_1'$. It follows that $\X_1$, $\U_1$ and $\U_1'$ have the same last model. Call it $\N$. It follows from full normalization that $\N$ is a $\Sigma^g$-iterate of $\S$ via a normal iteration that is above $\pi_{\Q_\b, \S}(\eta^\b)$. 

Because $\pi^\Q_{\R, \infty}(y)=\pi^\Q_{\R', \infty}(y')$, we have that $\pi_{\R, \N}(y)=\pi_{\R', \N}(y')$. Hence, $\pi_{\W, \S}(y)=\pi_{\W', \S}(y')$ (this uses the fact that $\T_{\S, \N}$ is above $\pi^{\X_0}(\eta^\b)$).  We now have the following equalities.

\begin{align*}
\pi^\P_{\W, \infty}(y)&=\pi^\P_{\S, \infty}(\pi_{\W, \S}(y))\\
&=\pi^\P_{\S, \infty}(\pi_{\W', \S}(y'))\\
&=\pi^\P_{\W', \infty}(y').
\end{align*}
\end{proof}

To finish the proof of the lemma, we show that $j$ is onto and $\Sigma_1$-elementary. Fix $u\in \M_\infty(\P)|\k_\infty^\P$. We want to see that there is $x\in \M_\infty(\Q)|\k_\infty^\Q$ such that $j(x)=u$. To start with, let $\R\in \mathcal{F}^g_\P$ be such that for some $v\in \R$, $\pi^\P_{\R, \infty}(v)=u$. Let $\a$ be such that $\T_{\P, \R}$ is based on $\P|\eta_\a$ and $v\in \R|\pi_{\P, \R}(\eta_\a)$. Using the fact that $\Sigma^g_{\P|\eta_\a}$ is directed, we can find an $\S$ that is a common $\Sigma^g_\R$ and $\Sigma_{\Q_\a}^g$ iterate of $\R$ and $\Q_\a$ respectively\footnote{A similar argument was used in the proof of Claim 2 to conclude that $\X_0$, $\U_0$ and $\U_0'$ have the same last model.}. Notice that $\T_{\Q_\a, \S}$ is based on $\Q_\a|\eta^\a$ and $\T_{\R, \S}$ is based on $\R|\pi_{\P, \R}(\eta_\a)$. Set $\X=(\T_{\Q_\a, \S}\rest (\Q_\a|\eta^\a))^\Q$. Letting $\N$ be the last model of $\X$, we  have that $\N$ is a normal strongly non-dropping $\Sigma^g_\S$-iterate of $\S$ via an iteration that is above $\pi_{\Q_\a, \S}(\eta^\a)$. Let $y=\pi_{\R, \S}(v)$. Notice that $\pi_{\S, \N}(y)=y$.  

We now easily have that $j(\pi^\Q_{\N, \infty}(y))=u$. Indeed, 
\begin{align*}
j(\pi^\Q_{\N, \infty}(y))&=\pi^\P_{\S, \infty}(y)\\&=\pi^\P_{\S, \infty}(\pi_{\R, \S}(v))\\ &=\pi^\P_{\R, \infty}(v)\\ &=u.
\end{align*}
To see that $j$ is $\Sigma_1$-elementary, it is enough to notice that given a $\Sigma_1$-formula $\phi(...)$, $\M_\infty(\Q)|\k_\infty^\Q\models \phi[\vec{x}]$ if and only if there is $\R\in \mathcal{F}^g_\Q$ such that for some $\vec{y}\in \R|\k$ and some $\a>\a_{\R, \vec{y}}$, $\pi^\Q_{\R, \infty}(\vec{y})=\vec{x}$ and 
\begin{center}
$\W(\R, \vec{y}, \a)|\pi_{\Q_\a, \W(\R, \vec{y}, \a)}(\eta^\a)\models \phi[\vec{y}]$. 
\end{center}
As $j(\vec{x})=\pi^\P_{\W(\R, \vec{y}, \a), \infty}(\vec{y})$, the $\Sigma_1$-elementarity follows. 

Since $j$ is onto and $\Sigma_1$-elementary, it follows that 
\begin{center}
$\M_\infty(\Q)|\k_\infty^\Q=\M_\infty(\P)|\k_\infty^\P$.
\end{center}
The fact that $\M_\infty(\Q)=\M_\infty(\P)$ follows from full normalization.
\end{proof}

Set $\k_\infty=\k^\P_\infty$ and $\M_\infty=\M_\infty(\P)$. The following useful corollary is an immediate consequence of excellence  and \rthm{1}. 

\begin{corollary}\label{useful corollary} Suppose $\R$ is a genericity iterate of $\P$ and $\Q$ is a genericity iterate of $\R$. Then $\Q$ is a genericity iterate of $\P$ (see \rprop{comute genericity}) and 
\begin{center}
$\pi^\R_{\R, \infty}=\pi^\Q_{\Q, \infty}\circ \pi_{\R, \Q}$.
\end{center}
\end{corollary}

The following proposition will be used in the next section. 

\begin{proposition}\label{2} Suppose $\Q$ is a genericity iterate of $\P$ and $h$ is a maximal $\Q$-generic. Suppose $X\in \V[g]$ is a countable subset of $\M_\infty(\P)|\xi$ for some $\xi<\k^\P_\infty$. Then $X\in \Q(\bR^*_h)$. 
\end{proposition}
\begin{proof} Let $(w_\a=(\eta_\a, \d_\a): \a<\cf(\k))$ and $(\Q_\a: \a<\cf(\k))$ witness that $\Q$ is a genericity iterate of $\P$. We can fix $\R\in \mathcal{F}^g_\P$, $\a<\cf(\k)$, a real $z\in \bR_g$ that codes a bijection $f:\omega\rightarrow \R|\pi_{\P, \R}(\eta_\a)$ and a $y\subseteq \omega$ such that 
\begin{center}
$\b\in X$ if and only if there is $i\in y$ such that $\pi^\P_{\R, \infty}(f(i))=\b$. 
\end{center}
Let $\S$ be a common iterate of $\Q_\a$ and $\R$. We have that $\S\in \mathcal{F}_{\Q_\a}^g$. Let $\N$ be the last model of  $(\T_{\Q_\a, \S}\rest \Q_\a|\eta^\a)^\Q$. Let $\sigma: \R|\pi_{\P, \R}(\eta_\a)\rightarrow \N|\pi_{\Q, \N}(\eta^\a)$ be the iteration embedding. We now have that 
\begin{center}
$\b\in X$ if and only if there is $i\in y$ such that $\pi^\Q_{\N, \infty}(\sigma(f(i)))=\b$. 
\end{center}
It follows from \rlem{-1} that $X\in \Q[z]$.
\end{proof}

\section{A Chang model over the derived model}

Suppose $\V\models \sf{ZFC}$ is a hod premouse and  $\k$ is a limit of Woodin cardinals of $\V$ such that if $\cf(\k)<\k$ then $\cf(\k)$ is not a measurable cardinal. Let $g\subseteq Coll(\omega, <\k)$ be $\V$-generic and let $\P=\V|(\k^+)^\V$. Recall that for $\b<\k$, $g_\b=g\cap Coll(\omega, <\b)$. Working in $\V(\bR^*_g)$, set
\begin{center}
$\bR^*_g=\cup_{\b<\k}\bR^{\V[g_\b]}$,\\
$\Gamma^*_g=\{A^g\cap \bR^*: \exists \b<\k (A\in \Gamma^\k_{g_\b})\}$,\\
$C(g)=L(\M_\infty, \cup_{\xi<\k_\infty}\powerset_{\omega_1}(\M_\infty|\xi), \Gamma^*_g, \bR^*_g)$,\\
$C_\a(g)=L_\a(\M_\infty, \cup_{\xi<\k_\infty}\powerset_{\omega_1}(\M_\infty|\xi), \Gamma^*_g, \bR^*_g)$
\end{center}

If $\Q\in \mathcal{I}^g(\P)$ then we let $\V_\Q$ be the last model of $(\T_{\P, \Q})^\V$. 

\begin{theorem}\label{determinacy in chang model} Suppose $(\V, \Omega)$ is an excellent pair, $\k$ is a limit of Woodin cardinals of $\V$ such that either $\k$ is regular or its cofinality is not a measurable cardinal, and $g\subseteq Coll(\omega, <\k)$ is $\V$-generic. Then $C(g)\models \sf{AD^+}$.
\end{theorem}

As with \rsec{internal dir lim cons: sec}, the use of $\Omega$ is to make sure that $S^\V$ has the desired properties. We will drop $\Omega$ from now on. 

\begin{proposition}\label{invariance} Suppose $\Q$ is a genericity iterate of $\P$ and $h$ is a maximal $\Q$-generic. Then $C(g)=(C(h))^{\V_\Q[h]}$.
\end{proposition}
\begin{proof} Because of \rthm{1} and \rprop{2}, it is enough to show that $\Gamma^*_g=(\Gamma^*_h)^{\V_\Q[h]}$. This is a standard fact about hod mice, and reduces to the fact that $\Sigma^g$ is pullback consistent. Indeed, suppose $B\in \Gamma^*_g$. Let $\eta<\k$ be an inaccessible cardinal such that the if $\d$ is the least Woodin cardinal of $\V$ above $\eta$ then there is a pair  of trees $(T, S)\in \V[g_\eta]$ witnessing that $B$ is $\k$-uB. Using genericity iterations, it can be shown that $B$ is projective in the code of $\Sigma^g_{\P|\d}$.  Given this information, notice next that $\Sigma^g_{\P|\d}\in (\Gamma^*_h)^{\V_\Q[h]}$ as it can be computed from $\Sigma^g_{\Q|\pi_{\P, \Q}(\d)}$ (using pullback consistency). This shows that $\Gamma^*_g\subseteq (\Gamma^*_h)^{\V_\Q[h]}$. That $(\Gamma^*_h)^{\V_\Q[h]} \subseteq \Gamma^*_g$ can be easily verified using the fact that $\Q$ is an iterate of $\P$. We now show that $B$ is indeed projective in $\Sigma^g_{\P|\d}$\footnote{That this is indeed the case is a standard fact and appears in many places in literature.}. Set $w=(\eta, \d)$.

Working in $\P[g_\eta]$, let $\tau\in \P[g_\eta]^{Coll(\omega, \d)}$ be the standard term relation that is always realized as $p[T]$ in $\P[g_\eta]^{Coll(\omega, \d)}$. 

We now have that $x\in B$ if and only if whenever $\S$ is a $\Sigma$-iterate of $\P$ based on $w$ such that $\pi_{\P, \S}$ is defined and $h\subseteq Coll(\omega, \pi_{\P, \S}(\d))$ is $\S[g_\eta]$-generic such that $x\in \S[g_\eta][h]$, $x\in (\pi_{\P, \S}(\tau))_h$. 

Indeed, let $x\in B$ and fix $(\S, h)$ as in the right side of the equivalence. Then since $x\in p[T]$, $x\in p[\pi_{\P, \S}(T)]$ and therefore, $x\in (\pi_{\P, \S}(\tau))_h$. Conversely, if $(\S, h)$ is as in the right side of the equivalence and $x\in (\pi_{\P, \S}(\tau))_h$, then $x\in p[\pi_{\P, \S}(T)]$. But then $x\in B$ as otherwise $x\in p[S]$, implying that $x\in  p[\pi_{\P, \S}(S)]$, contradiction. 
\end{proof}

We spend the rest of this section proving \rthm{determinacy in chang model}. To start with, fix $A\subseteq \bR^*_g$ and let $\a$ be the least such that $A\in C_\a( g)$. A Skolem hull argument shows that $\a<\k^{++}$\footnote{Indeed, working in $\V$, we can find $\b\geq \kappa^{+3}$ such that $A\in C_\b(g)$, and $\pi: M\rightarrow \V|\b$ such that $M^\omega\subseteq M$, $\card{M}=\kappa^+$ and $A\in M(\bR^*_g)$. But $C(g)^M=C_\a(g)$ where $\a=Ord\cap M$.}.  Let $(\gg, Y, B, x, s, \phi)$  be such that 
\begin{enumerate}
\item $\gg<\k_\infty$,
\item $Y\in \powerset_{\omega_1}(\M_\infty|\gg)$,
\item $B\in \Gamma^*_g$,
\item $x\in \bR^*_g$,
\item $s\in Ord^{<\omega}$,
\item for all $u\in \bR^*_g$,
\begin{center}
$u\in A \iff C_\a(g)\models \phi[Y, B, x, s, u]$.
\end{center}
\end{enumerate}
It is enough to show that $A\in \Gamma^*_g$, and for this, it is enough to show that for some $\l<\k$, $A$ is projective in $\Sigma^g_{\P|\l}$. Notice that there is some $\l_0$ such that $B$ is projective in $\Sigma_{\P|\l_0}^g$ (for instance, see the proof of \rprop{invariance}). Thus, without loss of generality we can assume that $B$ is the code of $\Sigma_{\P|\l_0}^g$ (denoted by $Code(\Sigma_{\P|\l_0}^g)$), and the real involved in the definition of our old $B$ is now part of $x$. We thus have that for every $u\in \bR^*_g$,\\\\
(1) $u\in A$ if and only if $C_\a(g)\models \phi[Y, Code(\Sigma_{\P|\l_0}^g), x, s, u]$.\\

Suppose $\R$ is a genericity iterate of $\P$. We say $\R$ is $A$-\textit{stable} if whenever $\Q$ is a genericity iterate of $\R$,  
\begin{center}
$\pi_{\R, \Q}[Y]=Y$ and $\pi_{\R, \Q}(\a, s)=(\a, s)$. 
\end{center}

\begin{lemma}
There is an $A$-stable $\R$.
\end{lemma}
\begin{proof}
First find a genericity iterate $\R_0$ of $\P$ such that $Y \subseteq rng(\pi^{\R_0}_{\R_0, \infty}[\k])$\footnote{This is one of the places where we use the fact that if $\k$ is singular then its cofinality is not a measurable cardinal. Otherwise there may not be such an $\R_0$ as $\pi_{\P, \infty}(\k)>\sup(\pi_{\P, \infty}[\k])$.}. Next find a genericity iterate $\R$ of $\R_0$ such that for any genericity iterate $\Q$ of $\R$, $\pi_{\R, \Q}(\a, s)=(\a, s)$. Notice that it follows from  \rcor{useful corollary} that we still have that 
$Y\subseteq rng(\pi^{\R}_{\R, \infty}[\k])$.

We have that $\R$ is a genericity iterate of $\P$. Moreover, whenever $\Q$ is a genericity iterate of $\R$, $\pi_{\R, \Q}[Y]=Y$. Indeed, fix such a $\Q$ and let $\xi\in \pi_{\R, \Q}[Y]$. Fix $\b\in Y$ such that $\xi=\pi_{\R, \Q}(\b)$ and let $\theta$ be such that $\pi_{\R, \infty}^\R(\theta)=\b$. But then 
\begin{align*}
\pi_{\R, \Q}(\xi) &=\pi_{\R, \Q}(\pi_{\R, \infty}^\R(\theta)) \\
&=\pi_{\Q, \infty}^\Q(\pi_{\R, \Q}(\theta))\ \text{(uses the fact that $\pi_{\R, \infty}^\R\in \V_\R$)}\\
&=\pi_{\R, \infty}^\R(\theta)\ (\text{see}\ \rcor{useful corollary})\\
&=\b.
\end{align*} Hence, $\xi=\b$ implying that $\xi\in Y$. The proof that $Y\subseteq  \pi_{\R, \Q}[Y]$ is similar. 
 \end{proof}
 We now fix some objects.
 \begin{enumerate}
 \item Let $\R$ be $A$-stable, and let $\tau<\k$ be such that $Y\subseteq \pi^\R_{\R, \infty}[\R|\tau]$. 
\item Let $u\in \bR^*_g$ be a real that codes $(\P|\l_0, \pi_{\P, \R}\rest (\P|\l_0), \R|\pi_{\P, \R}(\l_0))$. Notice that $\Sigma_{\P|\l_0}^g$ can be easily defined over $\R(\bR^*_g)$ from $u$ and $\Sigma_{\R|\pi_{\P, \R}(\l_0)}^g$ as the $\pi_{\P, \R}\rest (\P|\l_0)$-pullback of $\Sigma_{\R|\pi_{\P, \R}(\l_0)}^g$. 
\item Next, let $v\in \bR^*_g$ be a real coding a pair of reals $(v_0, v_1)$ such that $v_1$ codes a bijection $k:\omega \rightarrow \R|\tau$ and $Y=\{ \pi_{\R, \infty}^\R(k(n)): n\in v_0\}$.
\end{enumerate}
Notice that $Y\in \R[u, v]$ and $\R[u, v]\models ``Y$ is countable".

We can now get a formula $\psi$ such that for every $z\in \bR^*_g$, $z\in A$ if and only if 
\begin{center}
$\V_\R[x, u, v, z]\models \psi[\R|\k^{++}, \k, \a, s, x, u, v, z]$. 
\end{center}
The formula $\psi(W, p, q, r, a, b, c, d)$ essentially says the following (in the language of lbr hod mice):
\begin{enumerate}
\item $W$ is an lbr hod mouse,
\item $p$ is a limit of Woodin cardinals and $p^{++}=Ord\cap W$,
\item $q<p^{++}$ is an ordinal and $r$ is a finite sequence of ordinals,
\item $(a, b, c, d)\in \bR^4$,
\item $H_{p^{++}}$ is the universe of $W$,
\item $b$ codes a triple $(\M, i, \N)$ such that $i:\M\rightarrow \N$ is an elementary embedding and $\N\insegeq W|p$,
\item $c$ codes a pair of reals $(c_0, c_1)$ such that for some $\K\insegeq W|p$, $c_0$ codes an enumeration of $k:\omega\rightarrow \K$,
\item setting $Z=\{\pi^{W|p^+}_{W|p^+, \infty}(k(n)): n\in c_1\}$ and letting $\Lambda$ be the $i$-pullback of $\S^W_\N$, whenever $h\subseteq Coll(\omega, <p)$ is $W$-generic $(C_q(h))^{W(\bR^*_h)}\models \phi[Z, Code(\Lambda^h), a, r, d]$.
\end{enumerate} 

Let now $w$ be a window of $\R$ such that $\max(\pi_{\P, \R}(\l_0), \tau)<\eta^w$ and $(x, u, v)$ are generic for a poset in $\R|\eta^w$. We claim that $\l=\d^w$ is as desired, i.e., $A$ is projective in $\Sigma_{\R|\l}$. The following lemma shows exactly this.

Suppose $\S$ is an iterate of $\R$ via an iteration based on $w$ such that $\pi_{\R, \S}$ is defined. Let $\tau_\S\in \S[x, u, v]^{Coll(\omega, \l)}$ be the standard term  that is always realized as the set of reals $z$ such that 
\begin{center}
$\V_\S[x, u, v, z]\models \psi[\S|\k^{++}, \k, \a, s, x, u, v, z]$. 
\end{center}
We have that $\tau_\S\in \S[x, u, v]|(\l^{++})^\S$.

\begin{lemma} Suppose $z\in \bR^*_g$. Then $z\in A$ if and only if for any $\S\in \mathcal{F}^g_\R$ such that $\T_{\R, \S}$ is based on $w$  and $z$ is generic over $\S$ for ${\sf{EA}}^\S_{\pi_{\R, \S}(w)}$, and for any $\S[x, u, v]$-generic $h\subseteq Coll(\omega, \pi_{\R, \S}(\l))$ with the property that $z\in \S[x, u, v][h]$, $z\in (\tau_\S)_h$.\end{lemma}
\begin{proof} 
Assume $z\in A$ and let $(\S, h)$ be such that
\begin{enumerate}
\item $\S\in \mathcal{F}^g_\R$, $\T_{\R, \S}$ is based on $w$ and $z$ is generic over $\S$ for ${\sf{EA}}^\S_{\pi_{\R, \S}(w)}$, and
\item $h\subseteq Coll(\omega, \pi_{\R, \S}(\l))$ is $\S[x, u, v]$-generic with the property that $z\in \S[x, u, v][h]$.
\end{enumerate}
We need to see that $z\in (\tau_\S)_h$.  Continue $\T_{\R, \S}$ normally to get a $\Q$ that is a genericity iterate of $\R$ (see \rprop{constructing genericity iterations}). Because $\cp(\pi_{\R, \Q})>\eta^w$, we have that 
\begin{enumerate}
\item $u$ codes $(\P|\l_0, \pi_{\P, \Q}\rest \P|\l_0, \Q|\pi_{\P, \Q}(\l_0))$ and 
\item $v$ codes a pair of reals $(v_0, v_1)$ such that $v_1$ codes a bijection $k:\omega \rightarrow \Q|\tau$ and $Y=\{ \pi_{\Q, \infty}^\Q(k(n)): n\in v_0\}$.
\end{enumerate}

Notice now that because $z\in A$, we have that  $C_\a(g)\models \phi[Y, Code(\Sigma_{\P|\l_0}^g), x, s, z]$. Tracing definitions, we see that 
\begin{center}
$\V_\Q[x, u, v, z]\models \psi[\Q|\k^{++}, \k, \a, s, x, u, v, z]$. 
\end{center}
It then follows that $z\in (\tau_\Q)_h$. 
Because $\R$ is $A$-stable, we have that $\pi_{\R, \Q}(\a, s)=(\a, s)$ and $\pi_{\R, \Q}(Y)=\pi_{\R, \Q}[Y]=Y$\footnote{$\pi_{\R, \Q}(Y)$ makes sense as $Y\in \R[x, u, v]$ and $\cp(\pi_{\R, \Q})>\eta^w$.}, implying that $\pi_{\S, \Q}(\tau_\S)=\tau_\Q$. Because $\cp(\pi_{\S, \Q})>\pi_{\R, \S}(\l)$, it follows that $\tau_\S=\pi_{\S, \Q}(\tau_\S)$. Hence, $\tau_S=\tau_\Q$ and therefore, $x\in (\tau_\S)_h$.

Assume now that $z$ has the property that $\S\in \mathcal{F}^g_\R$ such that $\T_{\R, \S}$ is based on $w$  and $z$ is generic over $\S$ for ${\sf{EA}}^\S_{\pi_{\R, \S}(w)}$, and for any $\S[x, u, v]$-generic $h\subseteq Coll(\omega, \pi_{\R, \S}(\l))$ with the property that $z\in \S[x, u, v][h]$, $z\in (\tau_\S)_h$. We now want to see that $z\in A$.

Once again continue $\T_{\R, \S}$ normally to obtain a $\Q$ that is a genericity iterate of $\R$. We now have that $z\in (\tau_\Q)_h$ (because $\tau_\S=\pi_{\S, \Q}(\tau_\S)=\tau_\Q$). It then follows that
\begin{center}
$\Q[x, u, v, z]\models \psi[\Q|\k^{++}, \k, \a, s, x, u, v, z]$. 
\end{center}
Reversing the above implications, we get that $C_\a(g)\models \phi[Y, Code(\Sigma_{\P|\l_0}^g), x, s, z]$. Thus, $z\in A$.
\end{proof}

\section{The proof of \rthm{main theorem}}

Let $\M=\M_\infty$. We then have that $\nu=\k_\infty$. Clauses 1-3 of \rthm{main theorem} follow from \rthm{determinacy in chang model}.  Clause 4 follows from the results of \cite{TrangSteel}.

 \rcon{main conjecture} will eventually be settled by a robust core model induction technique which at the moment of writing this paper humanity does not posses.  However, it is an interesting question whether one can prove more instances of it in \textit{nice} structures. We suspect that it can be derived from just assuming that the universe has a unique universally Baire iteration strategy. For more on this concept see \cite{SargTrang}.

\bibliographystyle{plain}
\bibliography{CovuBrev.bib}
\end{document}